\documentclass{article}

\usepackage[utf8]{inputenc}
\usepackage{graphicx}
\usepackage{amsmath}
\usepackage{amssymb}
\usepackage{enumitem}
\usepackage{color}
\usepackage{tikz}
\usepackage{cases}
\usepackage{hyperref}
\usepackage{comment}
\usepackage{amsthm}
\usepackage{dsfont} 

\newtheorem{remark}{Remark}
\newtheorem{theorem}{Theorem}
\newtheorem{lemma}{Lemma}
\newtheorem{corollary}{Corollary}
\newtheorem{definition}{Definition}
\newtheorem{proposition}{Proposition}

\def\eps{\varepsilon}

\def\b0{{\bf 0}}
\def\e{{\bf e}}
\def\f{{\bf f}}
\def\g{{\bf g}}

\def\u{{\bf u}}
\def\v{{\bf v}}
\def\bvarphi{\boldsymbol{\varphi}}
\def\bpsi{\boldsymbol{\psi}}
\def\btau{\boldsymbol{\tau}}
\def\bsigma{\boldsymbol{\sigma}}

\def\A{\mathrm{\bf A}}
\def\B{\mathrm{\bf B}}
\def\D{\mathrm{\bf D}}
\def\S{\mathrm{\bf S}}

\def\Id{\mathrm{\bf I_3}}

\def\div{\mathrm{div}}

\title{Weak solution for granular model}

\author{
Laurent Chupin
\thanks{Universit\'e Clermont Auvergne, CNRS, LMBP, F-63000 Clermont-Ferrand, France\\
\phantom{toto}Corresponding author: laurent.chupin@uca.fr}  ~\& Thierry Dubois~$^*$
}

\date{}

\begin{document}

\setlength{\parindent}{0pt}

\maketitle


\begin{abstract}
This article is devoted to questions concerning the existence of solutions for partial differential equation problems modeling granular flows.
The models studied take into account the complex threshold rheology of these flows, as well as the dilatance effects. It is the coupling of these two physical phenomena that ensures stability and the existence of dissipated energy.
The key point of the article is to understand how this energy can ensure the existence of a weak solution.
We first establish a complete result on a simplified model, then demonstrate how it can be extended to more general cases.
This work represents a real breakthrough in the mathematical analysis of this type of models for complex flows.\\[0.3cm]
{\bf keywords}: granular model, weak existence, rheology, dilatance.
\end{abstract}

\section{Introduction}\label{sec:introduction}

\paragraph{From grain to continuous model}

Granular flow modeling is a major challenge in many contexts.
Indeed, this type of flow is involved in a wide variety of fields, such as the earth sciences, with pyroclastic flows or dune movements; the food industry in the broadest sense of the term, from cereal silos to food design; medicine, and in particular the pharmaceutical industry, which makes extensive use of powders and capsules corresponding to granular media with diverse behaviors.
The tricky aspect in this type of medium lies in its intermediate position between a fluid-like medium (such as gases, where the particles are too small to be individually described) and a solid medium, which requires understanding the interactions between a limited number of macroscopic particles.\\[0.2cm]
With the advent of computational tools, we could be led to use these models even when the media comprise a very large number of particles (generally several billion).
Currently this approach is hampered by our limits knowledge of the
local interactions: at the grain level, the solid contact laws between multiple particles involve highly non-linear relations such as frictional or inelastic shocks, that remain difficult to fully understand.
The well-known Discrete Element Method (see for instance~\cite{radjai11}) attempts to take into account as much information as possible, but the lack of physical understanding means that not all the expected results can be obtained.\\[0.2cm]
An alternative strategy is to define spatially averaged quantities to obtain continuous models, similar to what is done for fluid flows with the Navier-Stokes equations.
One of the foundational works of this approach for granular flows is arguably that of the GdR MIDI, which highlighted the $\mu(I)$-rheology model as a basis for granular flow models (see, for instance, the collective work~\cite{gdrmidi}).

\paragraph{Issues and open questions}

Since the $\mu(I)$-model is primarily derived from experimental considerations, numerous questions remain regarding its applicability in more general contexts those in which it was established (for instance, at larger scales).
Furthermore, its implementation raises additional challenges to be used numerically.
Indeed, to be able to numerically discretise a model and perform a computational implementation, it is highly desirable that it has a unique solution and that this solution is stable, in particular not too sensitive to perturbations in the data.
Unfortunately, from this point of view, these continuous models often exhibit significant limitations including issues of instability (see the work of Barker and colleagues in~\cite{barker23, barker17, barker15}) as well as the absence of a comprehensive theory ensuring the existence of solutions.
At first glance, these challenges appear to be closely linked to the non-linear and singular nature of the model.
Despite these fundamental issues, the implementation of these models in a numerical framework is feasible, but their reliability cannot be guaranteed. Indeed, instabilities have been frequently observed with classical formulations as reported in~\cite{chupin21}.
Several authors have presented results based on these models, but they have introduced assumptions that are not always physically relevant. This is the case, for example, in~\cite{abbatiello21,abbatiello19,chupin17} where the authors examine the existence of a solution in a model with a pressure-dependent stress threshold. However, the pressure considered in these models differs from the total pressure which results in a simpler coupling.
To address this limitation, one would need to incorporate more realistic models commonly used to describe granular flows (models capable of capturing phenomena such as normal stress differences, etc.). Such models have been developed and discussed in the review article~\cite{hutter94} and it is this type of model that we aim to explore in the present work.
\\[0.2cm]
The major challenge is therefore to develop a model for dense granular flows that is both physically consistent and mathematically well-posed. A first step in this direction was taken in~\cite{Chupin-Dubois24}, where a stable and physically consistent model was proposed: the particularly remarkable aspect is that it is the coupling between rheology and dilatancy that ensures the model's stability.
However, the strong non-linearities of the model, along with the presence of thresholds in the stress expression raise important questions regarding the existence and uniqueness of solutions.

\paragraph{Results and article outline}

In this article, we propose a natural functional framework in which it is possible to establish the existence and uniqueness of a solution. Once again, it is the combination of several physical ingredients that guarantees these results.
In the first part (Section~\ref{sec:model}), we describe the complete model and then focus on a simplified version that retains the two main physical features essential for ensuring the existence of a solution:
\begin{itemize}
\item[--] the incorporation of threshold rheology, with a pressure dependant yield criterion, as in the $\mu(I)$-model;
\item[--] the consideration of dilatancy, as proposed in~\cite{roux98}, in particular through a velocity field whose divergence, representing volume changes, is governed by a balanced between shear and pressure effects.
\end{itemize}
The second part form the core of the article and is devoted to the question of the existence of solutions for the proposed model. It is structured into several subsections.
After introducing the mathematical notations in Subsection~\ref{subsec:notations}, the following subsection presents a key component of the main result: it provides the weak formulation of the problem studied.
We observed in particular, that this formulation seems to eliminate certain threshold-related non-linearities. However, Proposition~\ref{prop:1} shows that this formulation effectively captures the full physical model and allows us to recover these non-linear effects.
The final subsections (Subsections~\ref{subsec:approximations},~\ref{subsec:estimates} and~\ref{subsec:limit}) are devoted to the proof of the existence result, with a particular focus on handling the remaining non-linear terms and establishing the positivity of pressure (a property that is physically intuitive, but rarely proven rigorously from a mathematical perspective).
These results represent a real innovation in the mathematical analysis of this type of model for complex flows: prove that a granular model, including threshold rheology and dilatation, is mathematically well-posed thanks to these relevant physical laws.\\[0.2cm]
In the final Section~\ref{sec:generalisation}, we return to the full problem and show how the preceding analysis can be adapted accordingly a specific study of the $\mu(I)$-model highlights how it relates to the “simplified” version initially considered.
This section concludes by pointing out that many questions remain open concerning this class of models.

\section{Modeling dense granular flow}\label{sec:model}

Guided by the following principles: make maximum use of physical principles (conservation laws) and phenomenological analyses (constitutive relations, equations of state); ensure the linear stability of the model, based on Barker's stability conditions presented in~\cite{schaeffer19}; and obtain an “energy” that dissipates over time, a complete dense granular flow model is proposed in~\cite{Chupin-Dubois24}.
This model couples four unknowns of the problem, namely the volume fraction scalar field~$\phi$, the velocity vector field~$\u$, the symmetric extra-stress tensor field~$\bsigma$ and the pressure scalar field~$p$.
Depending on the rheology chosen for the granular part, the expression of the various terms differs, but in the “simplest” case (corresponding to a Druker-Prager rheology), the model is written as follows
\begin{numcases}{}
\label{syst1} 
\phi \rho_0 (\partial_t \u + \u\cdot \nabla \u) + \nabla p - \div ( 2\nu(\phi, |\D\u|)\D\u) = \phi \rho_0 \g + \div \, \bsigma, \\
\label{syst3}
\bsigma:\S\u=2\alpha(\phi) p|\S\u|, \quad |\bsigma| \leq \alpha(\phi) p, \quad \bsigma=\bsigma^{\intercal} \quad\text{and}\quad \mathrm{tr}\, \bsigma=0, \\
\label{syst4}
\div \, \u = 2 \alpha(\phi)|\S\u| - \beta(\phi) \sqrt{p},\\
\label{syst0}
\partial_t \phi + \div(\phi\u) = 0.
\end{numcases}
In this version~$\g$ designates the gravity vector, $\rho_0$ is the grain density, $\phi_{\mathrm{max}}$ corresponds to the maximum volume fraction (of the order of~$0.6$ for experiment flows in laboratory, see~\cite{roche12}), $\alpha$ and~$\beta$ are positive functions depending on~$\phi$ and on several physical parameters such as average grain size, friction coefficient, etc.
The function~$\nu$ depends on the model choices: it reflects the possibly non-linear effect of viscous forces and depends on the strain-rate tensor corresponding to the symmetrical part of the velocity gradient:
$\D\u = \frac12 (\nabla \u + (\nabla \u)^\intercal)$.
The rheology also makes appear the deviatoric strain-rate tensor~$\S\u = \D\u - \frac13 (\div \, \u) \Id$. For all these tensors, the norm is defined by $|\A|^2 = \frac12 \A:\A$ where $\A: \B = \mathrm{tr}( \A\!^\intercal \, \B)$.\\[0.2cm]
The rheology introduced here corresponds to the Drucker-Prager rheology in the sense that the relations~\eqref{syst3} describing the relationship between stress and deformation can be rewritten as follows
\begin{equation}\label{syst-loi-seuil}\begin{aligned}
    & \bsigma = \alpha(\phi) p\frac{\S\u}{|\S\u|} \quad &&\text{if $|\S\u|\neq 0$},\\
    & |\bsigma|\leq \alpha(\phi) p, \quad \bsigma=\bsigma^{\intercal} \quad\text{and}\quad \mathrm{tr}\, \bsigma=0 \quad &&\text{if $|\S\u|= 0$}.
\end{aligned}\end{equation}
With this formulation, the quantity~$\alpha(\phi)$ can be seen as the threshold, often linked to the friction angle.\\[0.2cm] 
The other fundamental feature of the model~\eqref{syst1}--\eqref{syst0} is the expression for the divergence of the velocity field, i.e. the equation~\eqref{syst4}.
This is an interpretation of the law derived from the work of F. Radjaï and S. Roux in~\cite{roux98}. This law states that the local volume variation (i.e. $\div\, \u$) must be related to the volume fraction variation~$\phi$. More precisely, $\div\, \u$ must be proportional to the deviation between~$\phi$ and an equilibrium state~$\phi_{\mathrm{eq}}$ depending on~$|\S\u|$ and on~$p$, see Subsection~\ref{subsec:phi-variation} for more details.
\begin{remark}
As explained in~\cite{Chupin-Dubois24}, it is possible to choose more complex rheologies, such as the $\mu(I)$-rheology. In this case, however, the dilatation law~\eqref{syst4} must be adapted.
We will come back to the important case of $\mu(I)$-rheology in Section~\ref{sec:generalisation}.
However, for most granular models, the pressure~$p$ must be positive (for instance, the very definition of the inertial number~$I$ involves~$\sqrt{p}$, see Section~\ref{sec:generalisation}).
\end{remark}

\begin{remark}
This model does not take into account the interstitial gas in which the grains move.
In~\cite{Chupin-Dubois24}, this additional contribution is taken into account: essentially, it involves adding pressure forces~$\nabla p_f$ due to the gas in~\eqref{syst1}, and providing an evolution law for this new pressure~$p_f$. This evolution law is mainly a heat equation-type law, see~\cite{Chupin-Dubois24} for more details.
Note in particular that for this more complete model, as for the model presented here, energy is dissipated over time.
\end{remark}

\begin{remark}
In~\cite{Chupin-Dubois24}, the viscosity~$\nu$ is not taken into account, even though the authors mention this possibility in Remarks~3 and~4. In the case presented here, this viscosity is essential to ensure the regularity of the velocity field.
\end{remark}

In the following section, in order to highlight the tricky and original points, we focus on a simpler case of this model.
We will retain only the rheological and dilatational contributions, and we assume that the volume fraction is constant.
The model to be studied is the following
\begin{numcases}{}
\label{syst10} 
\phi_0 \rho_0 \partial_t \u + \nabla p - \div ( 2\nu_0|\D\u|\D\u) = \f + \div \, \bsigma, \\
\label{syst20}
\bsigma:\S\u=2\alpha_0 p|\S\u|, \quad |\bsigma| \leq \alpha_0 p, \quad \bsigma=\bsigma^{\intercal} \quad\text{and}\quad \mathrm{tr}\, \bsigma=0, \\
\label{syst30}
\div \, \u = 2 \alpha_0 |\S\u| - \beta_0 \sqrt{p},
\end{numcases}
where $\f$ are the external forces, namely~$\f=\phi_0 \rho_0 \g$.
Note that this model is not a special case of the full model, in which~$\phi$ would simply have been replaced by a constant~$\phi_0$, since the velocity field~$\u$ is not divergence free.
We show in Subsection~\ref{subsec:mu(I)} that this “simplified” system~\eqref{syst10}--\eqref{syst30} corresponds to the main order of the model~\eqref{syst1}--\eqref{syst0} in a specific regime, in particular when~$\phi$ is close to a constant itself close to~$\phi_{\mathrm{max}}$.
We propose in Section~\ref{sec:generalisation} an adaptation of these results to more complete cases, in particular taking into account the evolution of the volume fraction~$\phi$.
Finally, note the choice of non-linear viscosity. This choice is essentially technical, even if this type of viscosity may appear physically. In fact, as we shall see, it ensures sufficient velocity regularity to give meaning to the product~$p|\S\u|$.

\section{Mathematical results}\label{sec:maths}

\subsection{Mathematical framework and notations}\label{subsec:notations}

To carry out the mathematical study, we assume that all constants are normalized.
Of course, these assumptions do not qualitatively alter what follows.
The model considered in this section is therefore written as
\begin{numcases}{}
\label{pb1} 
\partial_t \u + \nabla p - \div ( 2|\D\u|\D\u) = \f + \div \, \bsigma, \\
\label{pb2}
\bsigma:\S\u=2p|\S\u|, \quad |\bsigma| \leq p, \quad \bsigma=\bsigma^{\intercal} \quad\text{and}\quad \mathrm{tr}\, \bsigma=0, \\
\label{pb3}
\div \, \u = 2 |\S\u| - \sqrt{p}.
\end{numcases}
Given a time $T>0$ and a bounded open domain $\Omega\subset \mathbb R^3$ whose boundary will be assumed to be of class~$\mathcal C^2$, we consider the problem~\eqref{pb1}--\eqref{pb3} on $(0,T)\times \Omega$.
It is completed by the following initial and boundary conditions:
\begin{equation}\label{pb4}
\u\big|_{t=0} = \u_{\mathrm{init}}
\qquad \text{and} \qquad
\u\big|_{\partial \Omega} = \b0.
\end{equation}

\paragraph{Notations -}
In order to define the notion of weak solution and to formulate the results, we need to fix the notations.
The symbol $\|\cdot\|_q$ for $1\leq q\leq +\infty$ stands for the $L^q$-norm in the usual Lebesgue space $L^q(\Omega)$ while $\|\cdot\|_{s,q}$ refers to the norm of the Sobolev space $W^{s,q}(\Omega)$, $s\in \mathbb R$ and $1\leq q\leq +\infty$.\\
We will make repeated use of the Sobolev space $W_0^{1,3}(\Omega)$ of functions of $W^{1,3}(\Omega)$ which are zero at the boundary, as well as its dual~$W^{-1,\frac32}(\Omega)$.
The associated duality bracket will be denoted $\langle\cdot,\cdot\rangle$.\\
Finally, for a Banach space $X$, we denote the relevant Bochner space by~$L^q(0,T;X)$. The associated norm will be explicitly denoted~$\|\cdot\|_{L^q(0,T;X)}$.

\begin{remark}|\label{rem:1756}
As explained in the previous section, the relation~\eqref{pb3} is derived from physical observations reflecting a competition between shear and pressure effects on volume changes.
Nevertheless, we can also see a fundamental mathematical ingredient that helps us understand the link between this relation~\eqref{pb3} and the rheology given by~\eqref{pb2}.
Indeed, like the role of pressure in the Stokes equations describing Newtonian fluid flow, the pressure~$p$ introduced in the equations~\eqref{pb1}--\eqref{pb2} seems to be linked to a Lagrange multiplier associated with the constraint relation~\eqref{pb3}.\\
As an example, consider the minimization problem $\u = \mathrm{argmin} \{f(\v)\,;\,g(\v)=\b0\}$
where
$$\begin{aligned}
f: \v \in W^{1,2}_0(\Omega) &\longmapsto && \frac{1}{2}\int_\Omega |\nabla \v|^2 \in \mathbb R\\
\qquad \text{and} \qquad
g: \v \in W^{1,2}_0(\Omega) &\longmapsto && \div \, \v - 2|\S\v| - q \in L^2(\Omega),
\end{aligned}
$$
where the source term~$q$ is given in $L^2(\Omega)$.
Note that the applications~$f$ and~$g$ are differentiable in any non-canceling function $\v\in W^{1,2}_0(\Omega)$ (if $\v$ cancels, we should probably consider the sub-differential of~$g$, as proposed by Beck~\cite{beck17}):
$$\begin{aligned}
df(\v): \bvarphi \in W^{1,2}_0(\Omega) &\longmapsto && \int_\Omega \nabla \v : \nabla \bvarphi \in \mathbb R\\
\qquad \text{and} \qquad
dg(\v): \bvarphi \in W^{1,2}_0(\Omega) &\longmapsto && \div \, \bvarphi - \frac{\S\v}{|\S\v|}:\S\bvarphi \in L^2(\Omega).\end{aligned}
    $$
If a solution~$\u$, minimizing~$f$ under the constraint $g(\u)=\b0$, exists and does not cancel, then according to the Lagrange multiplier Theorem, there exists a linear application $\widetilde{p}:L^2(\Omega) \longrightarrow \mathbb R$ such that $df(\u)=\widetilde{p}\circ dg(\u)$.
According to the Riesz representation Theorem, the linear application~$\widetilde{p}$ can be represented by the scalar product in~$L^2(\Omega)$: there exists $p\in L^2(\Omega)$ such that for all $\bpsi\in L^2(\Omega)$ we have $\widetilde{p}(\bpsi) = \int_\Omega p \bpsi$.
Thus, using the expressions for the differentials of the applications~$f$ and~$g$, we deduce that
$$
\forall \bvarphi\in W^{1,2}_0(\Omega) \quad \int_\Omega \nabla \u : \nabla \bvarphi = \int_\Omega p \big( \div \, \bvarphi - \frac{\S \u}{|\S\u|} : \S\bvarphi \big)
$$
which indicates that, in the sense of distributions, we have $-\Delta \u + \nabla p = \div\, \bsigma$ with $\bsigma:\S\u = 2p|\S\u|$.
\end{remark}

The approach described in Remark~\ref{rem:1756}, which allows pressure to be seen as a Lagrange multiplier, will not be used in this article, but could perhaps provide some interesting results.
We propose here an original formulation in terms of velocity-pressure, in which we emphasize the energy aspect.

\subsection{Weak formulation and statement of result}\label{subsec:weak}

An important point is the definition of the notion of solution of the system of equations~\eqref{pb1}--\eqref{pb4}.
\begin{definition}[weak solution]\label{definition}
Let $\u_{\mathrm{init}}\in L^2(\Omega)$ and $\f\in L^{\frac32}(0,T;W^{-1,\frac32}(\Omega))$.\\
We say that $(\u,p,\bsigma)$ is a weak solution of~\eqref{pb1}--\eqref{pb4} if
\begin{equation*}\begin{aligned}
& \u \in L^\infty(0,T;L^2(\Omega))\cap L^3(0,T;W_0^{1,3}(\Omega)), \quad \partial_t \u \in L^{\frac32}(0,T;W^{-1,\frac32}(\Omega)),\\
& p \in L^{\frac32}((0,T)\times \Omega) \quad \text{with $p\geq 0$ \, a.e.}\\
& \bsigma \in L^{\frac32}((0,T)\times \Omega) \quad \text{with}~ |\bsigma|\leq p, \quad \bsigma=\bsigma^{\intercal} \quad\text{and}\quad \mathrm{tr}\, \bsigma=0 ~~ \text{a.e.}
\end{aligned}\end{equation*}
and for all $\bvarphi\in L^3(0,T;W_0^{1,3}(\Omega))$, for all $\psi\in L^{\frac32}((0,T)\times \Omega)$ we have
\begin{align}
& \langle \partial_t \u , \bvarphi \rangle - \int_\Omega p\, \div \, \bvarphi + \int_\Omega 2|\D\u|\D\u:\D\bvarphi + \int_\Omega \bsigma:\S\bvarphi = \langle \f , \bvarphi \rangle, \label{weak1}\\
& \int_\Omega \psi \, \div \, \u - \int_\Omega 2\psi |\S\u| + \int_\Omega \psi \sqrt{p} = 0, \label{weak2}\\
& \frac{1}{2}\int_\Omega |\u|^2
+
4 \int_0^T\int_\Omega |\D\u|^3
+
\int_0^T\int_\Omega |p|^{\frac32}
\leq
\frac{1}{2}\int_\Omega |\u_{\mathrm{init}}|^2
+
\int_0^T \langle \f , \u \rangle,\label{weak3}
\end{align}
and the initial condition $\u\big|_{t=0} = \u_{\mathrm{init}}$ holds in $L^2(\Omega)$.
\end{definition}
Relationships~\eqref{weak1} and~\eqref{weak2} naturally translate equations~\eqref{pb1} and~\eqref{pb3} into their weak form. However it is not completely clear that threshold rheology as expressed in relations~\eqref{pb2} is fully taken into account in this notion of solution.
\begin{proposition}\label{prop:1}
If $(\u,p,\bsigma)$ is a weak solution of~\eqref{pb1}--\eqref{pb4} then the relations~\eqref{pb2} are satisfied almost everywhere.
\end{proposition}
{\bf Proof -}
Let $(\u,p,\bsigma)$ be a weak solution of~\eqref{pb1}--\eqref{pb4}.
Taking $\bvarphi=\u$ as test function in~\eqref{weak1}, $\psi=p$ in~\eqref{weak2} and adding the results, we deduce
\begin{equation*}
    \frac{1}{2}\frac{d}{dt}\int_\Omega |\u|^2 + \int_\Omega 4|\D\u|^3 + \int_\Omega |p|^{\frac32} + \int_\Omega \Big( \bsigma:\S\u - 2p|\S\u| \Big) = \langle \f , \u \rangle.
\end{equation*}
Integrating with respect to time and comparing with~\eqref{weak3}, we get
\begin{equation}\label{eq:2200}
    \int_0^T\int_\Omega \Big( \bsigma:\S\u - 2p|\S\u| \Big) \geq 0.
\end{equation}
But, from the Cauchy-Schwarz inequality $\bsigma:\S\u \leq 2|\bsigma| |\S\u|$ and using the fact that $|\bsigma|\leq p$, we have $\bsigma:\S\u - 2p|\S\u|\leq 0$. Regarding~\eqref{eq:2200} that implies $\bsigma:\S\u = 2p|\S\u|$ almost everywhere.
\hfill$\blacksquare$

\vspace{0.2cm}

The main result is 
\begin{theorem}\label{th:1}
If $\u_{\mathrm{init}}\in L^2(\Omega)$ and $\f\in L^{\frac32}(0,T;W^{-1,\frac32}(\Omega))$ then there exists a weak solution $(\u,p,\bsigma)$ to the problem~\eqref{pb1}--\eqref{pb4}.
The velocity~$\u$ and the pressure~$p$ are uniquely determined.
\end{theorem}

\subsection{Approximations system}\label{subsec:approximations}

In order to obtain a weak solution, we start by constructing a solution to an approximate problem. To do this, consider two sequences~$(\f_\eps)_{\eps>0}$ and~$(\u_{\mathrm{init},\eps})_{\eps>0}$ of regular functions respectively converging to~$\f$ in $L^\frac32(0,T;W^{-1,\frac32}(\Omega))$ and $\u_{\mathrm{init}}$ in $L^2(\Omega)$.
For $\eps>0$ we consider the following problem
\begin{numcases}{}
\label{pb1k}
\partial_t \u_{\eps} + \nabla p_{\eps} - \div ( 2|\D\u_{\eps}|\D\u_{\eps}) = \f_{\eps} + \div \, \bsigma_{\eps}, \\
\label{pb2k}
\bsigma_{\eps} = p_{\eps}\frac{\S\u_{\eps}}{|\S\u_{\eps}|+\eps}, \\
\label{pb3k}
\div \, \u_{\eps} = 2\frac{|\S\u_{\eps}|^2}{|\S\u_{\eps}|+\eps} - V_\eps(p_{\eps}) - \eps (\partial_t p_{\eps}-\Delta p_{\eps}).
\end{numcases}
with the following initial and boundary conditions:
\begin{equation}\label{pb4k}
\u_{\eps}\big|_{t=0} = \u_{\mathrm{init},\eps},
\qquad
\u_{\eps}\big|_{\partial \Omega} = \b0,
\qquad
p_{\eps}\big|_{t=0} = 0
\qquad \text{and} \qquad
p_{\eps}\big|_{\partial \Omega} = 0.
\end{equation}
The function~$V_\eps$ is an approximation of the square-root function, which allows us to define pressure~$p_\eps$ for all $\eps>0$, and to show that we obtain a positive pressure at the limit~$\eps\to 0$.
More precisely, $V_\eps$ is a concave function of class~$\mathcal C^1$, defined on~$\mathbb R$ by
\begin{equation}\label{defV}
V_\eps(x) = \left\{ \begin{aligned}
& \sqrt{x+(\eps/2)^2} - \eps/2 && \quad \text{if $x>0$},\\
& x/\eps && \quad \text{if $x\leq 0$}.
\end{aligned}\right.
\end{equation}
The following two graphs illustrate how the square root function and the Heaviside function have been approximated by smooth functions.
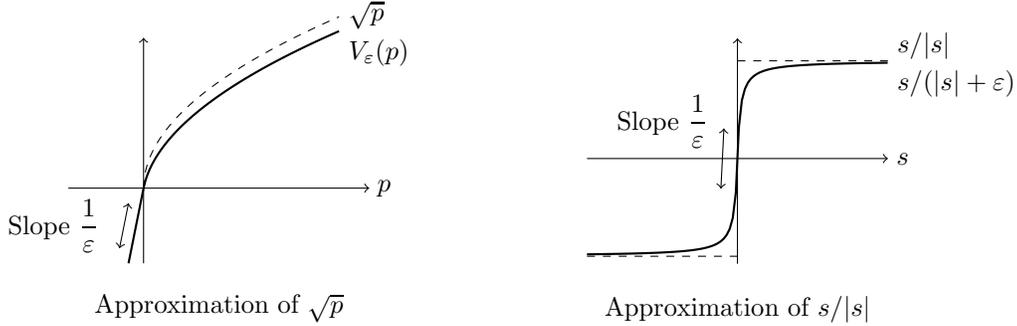
\begin{center}
\begin{tikzpicture}[scale=2]
\draw [->] (-0.5,0) -- (1.5,0) ;
\draw (1.6,0) node{$p$};
\draw [->] (0,-0.5) -- (0,1) ;
\draw [domain=0:1.3,samples=100, dashed] plot (\x,{sqrt(\x)});
\draw (1.3,1.15) node[right]{$\sqrt{p}$};
\draw [domain=0:1.3,samples=100, thick] plot (\x,{sqrt(\x+0.1*0.1)-0.1});
\draw (1.3,0.9) node[right]{$V_\eps(p)$};
\draw [thick] (-0.1,-0.5) -- (0,0) ;
\draw [<->] (-0.1,-0.1) -- (-0.1-0.3/5,-0.4) ;
\draw (-0.6,-0.25) node{Slope $\displaystyle \frac{1}{\eps}$};
\draw (0.5,-0.8) node{Approximation of $\sqrt{p}$};
\end{tikzpicture}
\hspace{2cm}
\begin{tikzpicture}[scale=2]
\draw [->] (-1,0) -- (1,0) ;
\draw (1.1,0) node{$s$};
\draw [->] (0,-0.7) -- (0,0.8) ;
\draw [dashed] (-1,-0.65) -- (0,-0.65);
\draw [dashed] (0,0.65) -- (1,0.65);
\draw (1,0.75) node[right]{$\displaystyle s/|s|$};
\draw [domain=-1:1,samples=100, thick] plot (\x,{0.65*\x/(abs(\x)+0.02)});
\draw (1,0.5) node[right]{$\displaystyle s/(|s|+\eps)$};
\draw [<->] (-0.1-0.01,-0.01*20) -- (-0.1+0.01,0.01*20) ;
\draw (-0.5,0.25) node{Slope $\displaystyle \frac{1}{\eps}$};
\draw (0,-1.) node{Approximation of $s/|s|$};
\end{tikzpicture}
\end{center}
\begin{proposition}\label{prop-existencek}
For $\eps>0$, there exists a unique weak solution $(\u_{\eps},p_{\eps},\bsigma_{\eps})$ to equations~\eqref{pb1k}--\eqref{pb4k}, that is such that
\begin{equation*}\begin{aligned}
& \u_{\eps} \in L^\infty(0,T;L^2(\Omega))\cap L^3(0,T;W_0^{1,3}(\Omega)), \quad &&\partial_t \u_{\eps} \in L^{\frac32}(0,T;W^{-1,\frac32}(\Omega)),\\
& p_{\eps} \in L^\infty(0,T;L^2(\Omega))\cap L^2(0,T;W_0^{1,2}(\Omega)), \quad &&\partial_t p_{\eps} \in L^2(0,T;W^{-1,2}(\Omega)),\\
& \bsigma_{\eps} \in L^{\frac32}((0,T)\times \Omega),
\end{aligned}\end{equation*}
with $\displaystyle \bsigma_{\eps} = p_{\eps}\frac{\S\u_{\eps}}{|\S\u_{\eps}|+\eps}$
and for all $\bvarphi\in L^3(0,T;W_0^{1,3}(\Omega))$, for all $\psi\in L^2(0,T;W_0^{1,2}(\Omega))$
\begin{align}
& \langle \partial_t \u_{\eps} , \bvarphi \rangle
-
\int_\Omega p_{\eps}\, \div \, \bvarphi
+
\int_\Omega 2|\D\u_{\eps}|\D\u_{\eps}:\D\bvarphi
+
\int_\Omega \bsigma_{\eps}:\S\bvarphi
=
\langle \f_{\eps} , \bvarphi \rangle, \label{weak1k}\\
& \langle \eps \partial_t p_\eps , \psi \rangle
+
\int_\Omega \eps \nabla p_\eps \cdot \nabla \psi
+
\int_\Omega \psi \, \div \, \u_{\eps}
-
\int_\Omega 2\psi \frac{|\S\u_{\eps}|^2}{|\S\u_{\eps}|+\eps}
+
\int_\Omega \psi V_\eps(p_{\eps}) = 0. \label{weak2k}
\end{align}
The initial conditions are satisfied in~$L^2(\Omega)$: $\u_{\eps}|_{t=0} = \u_{\mathrm{init},\eps}$ and $p_{\eps}|_{t=0} = 0$.
\end{proposition}

{\bf Proof -}
Adding $\partial_t p_{\eps}-\Delta p_{\eps}$ into the equation~\eqref{pb3k} allows us to define, for each $\eps>0$, the pressure~$p_{\eps}$ as the solution to a heat-type equation.
The explicit expression of $\sigma_{\eps}$ given by~\eqref{pb2k} is a usual regularized version of the condition~\eqref{pb2}. This type of regularization is very common in this kind of problem, see for example~\cite{abbatiello19}.
The existence of a unique solution to the complete system~\eqref{pb1k}--\eqref{pb3k} is then classical (for example, using a Galerkin method, see~\cite{abbatiello19}).
In practice, the Galerkin method applies as soon as we have an energy estimate for the solution.
We therefore show below how to derive such an estimate.\\[0.2cm]
We choose $\bvarphi=\u_{\eps}$ as test function in~\eqref{weak1k}, $\psi=p_{\eps}$ in~\eqref{weak2k}, we integrate with respect to time and add the results.
Note, in particular that, based on the expression for~$\bsigma_{\eps}$, we have
\[\bsigma_{\eps}:\S\u_{\eps} = 2p_{\eps} \frac{|\S\u_{\eps}|^2}{|\S\u_{\eps}|+\eps}.\]
We also note that, for $\u_{\eps} \in L^3(0,T;W_0^{1,3}(\Omega))$ and $\partial_t \u_{\eps} \in L^{\frac32}(0,T;W^{-1,\frac32}(\Omega))$, we know that (see~\cite[p.99]{boyer12})
\[\langle \partial_t \u_{\eps} , \u_{\eps} \rangle = \frac{1}{2} \frac{d}{dt} \|\u_{\eps}\|_2^2.\]
The same reasoning applied to the time derivative of the pressure term finally implies
\begin{align}
&
\frac{1}{2}\frac{d}{dt}\|\u_{\eps}\|_2^2
+
\frac{\eps}{2}\frac{d}{dt}\|p_{\eps}\|_2^2
+
4\|\D\u_{\eps}\|_3^3
+
\eps \|\nabla p_{\eps}\|_2^2
+
\int_\Omega p_\eps V_\eps(p_{\eps})
= 
\langle \f_{\eps} , \u_{\eps} \rangle.\label{weak3k}
\end{align}
This estimate allows us to demonstrate the existence of a solution to equations~\eqref{pb1k}--\eqref{pb4k} in the classical way. Note in particular that the choice of the function~$V_\eps$ ensures that for all $x\in \mathbb R$, we have $xV_\eps(x)\geq 0$.\\
The proof of uniqueness is also standard, and we refer the reader to Subsection~\ref{subsec:uniqueness} for a similar reasoning leading to the result.
\hfill$\blacksquare$

\subsection{Uniform estimates}\label{subsec:estimates}

To take the limit as $\eps\to 0$ in the system~\eqref{pb1k}--\eqref{pb4k}, we will establish uniform bounds with respect to~$\eps$ for the unknowns~$\u_\eps$, $p_\eps$ and $\bsigma_{\eps}$:
\begin{proposition}\label{prop-bounds}
The following quantities are bounded independently of $\eps>0$:
\begin{equation*}\begin{aligned}
    & \|\u_{\eps}\|_{L^\infty(0,T;L^2(\Omega))},\,
    &&\|\nabla \u_{\eps}\|_{L^3((0,T)\times \Omega)},\,
    &&\|\partial_t \u_{\eps}\|_{L^{\frac32}(0,T;W^{-1,\frac32}(\Omega))},\,
    &&\|\bsigma_{\eps}\|_{L^{\frac32}((0,T)\times \Omega)},\\
    &\|p_{\eps}\|_{L^{\frac32}((0,T)\times \Omega)},\,
    && \|\eps^{\frac{1}{2}} p_{\eps}\|_{L^\infty(0,T;L^2(\Omega))},\
    &&\|\eps^{\frac{1}{2}} \nabla p_{\eps}\|_{L^2((0,T)\times\Omega)},\,
    &&\|\eps \partial_t p_{\eps}\|_{L^2(0,T;W^{-1,2}(\Omega))}.
\end{aligned}\end{equation*}
\end{proposition}

{\bf Proof - Step 1: using energy estimate}\\[0.2cm]
We first estimate the term in the right-hand side of equation~\eqref{weak3k}.
Thus, using successively the duality between~$W^{-1,\frac{3}{2}}(\Omega)$ and~$W_0^{1,3}(\Omega)$, the Korn's inequality $\|\u\|_{1,3} \leq K_3 \|\D\u\|_3$, see \cite[Chapter V, Theorem 1.10]{malek96}, and the Young's inequality $ab \leq \frac{2}{9} a^{\frac32} + 3 b^3$, we obtain
\[\begin{aligned}
\langle \f_{\eps} , \u_{\eps} \rangle
& \leq \|\f_{\eps}\|_{-1,\frac32} \|\u_{\eps}\|_{1,3} \\
& \leq K_3 \|\f_{\eps}\|_{-1,\frac32}\|\D\u_{\eps}\|_{3} \\
& \leq \frac{2}{9}K_3^{\frac32} \|\f_{\eps}\|_{-1,\frac32}^{\frac32} + 3 \|\D\u_{\eps}\|_3^3.
\end{aligned}\]
Finally, the equality of energy~\eqref{weak3k} implies the following estimate:
\begin{equation*}\begin{aligned}
\frac{d}{dt}\|\u_{\eps}\|_2^2 
+
\eps \frac{d}{dt}\|p_{\eps}\|_2^2
+
2\|\D\u_{\eps}\|_3^3
+
2\eps \|\nabla p_{\eps}\|_2^2
+
2\int_\Omega p_\eps V_\eps(p_{\eps})
\leq
\frac{4}{9}K_3^{\frac32} \|\f_\eps\|_{-1,\frac{3}{2}}^{\frac{3}{2}}.
\end{aligned}\end{equation*}
Integrating in time and using the convergence of $\f_{\eps}$ to~$\f$ in $L^{\frac{3}{2}}(0,T;W^{-1,\frac{3}{2}}(\Omega))$, we deduce
\begin{equation}\label{estk-1}\begin{aligned}
& \sup_{(0,T)}\|\u_{\eps}\|_2^2
+
\eps \sup_{(0,T)} \|p_{\eps}\|_2^2
+
2 \int_0^T \|\D\u_{\eps}\|_3^3
+
2 \eps \int_0^T \|\nabla p_{\eps}\|_2^2
+
2 \int_0^T \int_\Omega p_\eps V_\eps(p_{\eps})\\
& \hspace{8cm}\leq
\|\u_{\mathrm{init}}\|_2^2 + K_3^{\frac32} \int_0^T  \|\f\|_{-1,\frac{3}{2}}^{\frac{3}{2}}.
\end{aligned}\end{equation}
Noting that, by definition of the function~$V_\eps$, see~\eqref{defV}, the quantity~$p_\eps V_\eps(p_{\eps})$ is always non-negative we therefore have estimates of the following quantities, independently of~$\eps$:
\begin{equation*}
    \|\u_{\eps}\|_{L^\infty(0,T;L^2(\Omega))}, \quad
    \|\eps^{\frac{1}{2}} p_{\eps}\|_{L^\infty(0,T;L^2(\Omega))}, \quad
    \|\D\u_{\eps}\|_{L^3((0,T)\times \Omega)} \quad \text{and} \quad
    \|\eps^{\frac{1}{2}} \nabla p_{\eps}\|_{L^2((0,T)\times\Omega)}.
\end{equation*}
In addition, the estimate of the quantity~$p_\eps V_\eps(p_{\eps})$ provides a bound on the pressure~$p_{\eps}$. Indeed, we have
\begin{equation*}
\int_0^T\int_\Omega |p_{\eps}|^{\frac{3}{2}} =
\underbrace{\iint_{0\leq p_{\eps}\leq 1} |p_{\eps}|^{\frac32}}_{A}
+
\underbrace{\iint_{p_{\eps}<0} |p_{\eps}|^{\frac32}}_{B}
+
\underbrace{\iint_{p_{\eps}>1} |p_{\eps}|^{\frac32}}_{C}.
\end{equation*}
We clearly have $A \leq T|\Omega|$.
Moreover, due to Hölder's inequality, we obtain
\begin{equation*}
B
\leq \Big(\iint_{p_{\eps}<0} \eps^3 \Big)^{\frac{1}{4}} \Big(\iint_{p_{\eps}<0} \frac{1}{\eps}|p_{\eps}|^2 \Big)^{\frac{3}{4}}
\leq (T|\Omega|\eps^3)^{\frac{1}{4}} \Big(\int_0^T \int_\Omega p_{\eps}V_\eps(p_\eps) \Big)^{\frac{3}{4}}.
\end{equation*}
Using the fact that, for $x\geq 1$ and $\eps\leq 2$ we have $x^{\frac{3}{2}} \leq \frac{1}{\sqrt{2}-1}xV_\eps(x)$, this implies that
\begin{equation*}
C \leq \frac{1}{\sqrt{2}-1} \int_0^T \int_\Omega p_{\eps}V_\eps(p_\eps).
\end{equation*}
These three estimates on~$A$, $B$ and~$C$ combined with the energy estimate~\eqref{estk-1} provide a bound on $\|p_{\eps}\|_{L^\frac{3}{2}((0,T)\times \Omega)}$ that is independent of~$\eps$, as long as $\eps\leq 2$).\\
The estimate on the stress~$\bsigma_{\eps}$ follows directly from that of the pressure, since expression~\eqref{pb2k} yields the inequality~$|\bsigma_{\eps}|\leq |p_{\eps}|$.\\[0.2cm]
{\bf Step 2: Bounds for times derivatives}\\[0.2cm]
Let us go back to the weak formulation~\eqref{weak1k}.
More precisely, using Hölder's inequality and the duality $W^{-1,\frac{3}{2}}(\Omega)\leftrightarrow W_0^{1,3}(\Omega)$, we know that for all $\bvarphi\in L^3(0,T;W_0^{1,3}(\Omega))$ we have
\begin{equation*}
    \langle \partial_t \u_{\eps} , \bvarphi \rangle
    \leq
    \|p_{\eps}\|_{\frac32} \|\div \, \bvarphi\|_3 + 2\|\D\u_{\eps}\|_3^2 \|\D\bvarphi\|_3 + \|\bsigma_{\eps}\|_{\frac32} \|\S\bvarphi\|_3 +  \|\f_{\eps}\|_{-1,\frac32} \|\bvarphi\|_{1,3}.
\end{equation*}
Integrating in time and applying Hölder's inequalities (with respect to integration in time), we derive the following estimate
\begin{equation*}\begin{aligned}
    \Big| \int_0^T \langle \partial_t \u_{\eps} , \bvarphi \rangle \Big|
    \leq
& \Big(\int_0^T \|p_{\eps}\|_{\frac32}^{\frac32} \Big)^{\frac23} \Big(\int_0^T \|\div \, \bvarphi\|_3^3 \Big)^{\frac13} 
    + 2\Big(\int_0^T \|\D\u_{\eps}\|_3^3\Big)^{\frac23} \Big(\int_0^T \|\D\bvarphi\|_3^3 \Big)^{\frac13} \\
& \qquad + \Big(\int_0^T \|\bsigma_{\eps}\|_{\frac32}^{\frac32} \Big)^{\frac23} \Big(\int_0^T \|\S\bvarphi\|_3^3 \Big)^{\frac13} 
    +  \Big(\int_0^T \|\f_{\eps}\|_{-1,\frac32}^\frac32 \Big)^{\frac23} \Big(\int_0^T \|\bvarphi\|_{1,3}^3 \Big)^{\frac13}.
\end{aligned}\end{equation*}
By virtue of the estimates~\eqref{estk-1}, we deduce that there exists a constant~$C_1$, independent of~$\eps$, such that
\begin{equation*}
    \Big| \int_0^T \langle \partial_t \u_{\eps} , \bvarphi \rangle \Big|
    \leq
    C_1 \Big(\int_0^T \|\bvarphi\|_{1,3}^3 \Big)^{\frac13} = C_1 \|\bvarphi\|_{L^3(0,T;W_0^{1,3}(\Omega))}.
\end{equation*}
Thus, by the definition of the dual norm, we obtain $\|\partial_t \u_{\eps}\|_{L^{\frac32}(0,T;W^{-1,\frac32}(\Omega))}\leq C_1$.\\[0.2cm]
Similarly, we use~\eqref{pb3k} to estimate $\partial_t p_{\eps}$. For all $\psi\in L^2(0,T;W_0^{1,2}(\Omega))$, we have
\begin{equation*}
\langle \eps \partial_t p_\eps , \psi \rangle
\leq
\eps \|\nabla p_\eps\|_2 \|\nabla \psi\|_2
+
\|\psi\|_\frac{3}{2} \|\div \, \u_{\eps}\|_3
+
2\|\psi\|_\frac{3}{2} \|\S\u_{\eps}\|_3
+
\|\psi\|_\frac{3}{2} \|V_\eps(p_{\eps})\|_3.
\end{equation*}
We integrate in time and use the previously established bounds on $\eps^{\frac{1}{2}}\nabla p_\eps$, $\div \, \u_\eps$ and $\S\u_\eps$ (the latter two relying on the bound for~$\D\u_\eps$).
It remains to estimate~$\|V_\eps(p_{\eps})\|_3$.\\
To do this, we take $\psi=|V_\eps(p_{\eps})|V_\eps(p_{\eps})$ in~\eqref{weak2k} and we obtain
\begin{align}
&
\frac{d}{dt}\int_\Omega \eps W_\eps(p_\eps)
+
\int_\Omega \eps W''_\eps(p_\eps) |\nabla p_\eps|^2 
+
\|V_\eps(p_{\eps})\|_3^3
\leq
\int_\Omega 2V_\eps(p_{\eps})^2 |\S\u_{\eps}|
+ \int_\Omega V_\eps(p_{\eps})^2 \, |\div \, \u_{\eps}|,
\label{weak3k:0828}
\end{align}
where $W_\eps$ is such that $W_\eps'=|V_\eps|V_\eps$ and $W_\eps(0)=0$.
The right-hand side is bounded using the previous estimates.
For instance, applying successively Hölder's inequality followed by Young's inequality, we derive
$$
\int_\Omega V_\eps(p_{\eps})^2 \, |\div \, \u_{\eps}|
\leq
\|V_\eps(p_{\eps})\|_3^2\|\div \, \u_{\eps}\|_3
\leq
\frac{1}{3}\|V_\eps(p_{\eps})\|_3^3 + \frac{2}{3}\|\div \, \u_{\eps}\|_3^3.
$$
Since $W_\eps$ and $W''_\eps$ are non-negative, and $W_\eps(p_\eps|_{t=0})=W_\eps(0)=0$, integration of~\eqref{weak3k:0828} with respect to time over the interval~$(0,T)$ yields the desired bound on~$V_\eps(p_\eps)$.\\
We deduce that there exists a constant~$C_2$, independent of~$\eps$, such that $\|\eps \partial_t p_{\eps}\|_{L^2(0,T;W^{-1,2}(\Omega))}\leq C_2$.
\hfill$\blacksquare$

\subsection{Limit process}\label{subsec:limit}

By virtue of the bounds obtained in Proposition~\ref{prop-bounds}, we deduce the following convergences (up to sub-sequence extractions):
\begin{proposition}\label{prop:cv}
There exists a function~$\u$ such that
\begin{equation*}\begin{aligned}
    & \u_{\eps} \rightharpoonup \u \quad &&\text{weakly-$\star$ in $L^\infty(0,T;L^2(\Omega))$},\\
    & \nabla\u_{\eps} \rightharpoonup \nabla\u \quad &&\text{weakly in $L^3(0,T;L^3(\Omega))$},\\
    & \partial_t \u_{\eps} \rightharpoonup \partial_t \u \quad &&\text{weakly in $L^{\frac32}(0,T;W^{-1,\frac32}(\Omega))$},\\
    & \u_{\eps} \rightarrow \u \quad &&\text{in $L^3(0,T;L^2(\Omega))$}.
\end{aligned}\end{equation*}
There exists a non-negative function~$p$ such that
\begin{equation*}\begin{aligned}
    & p_{\eps} \rightharpoonup p \quad &&\text{weakly in $L^{\frac32}(0,T;L^{\frac32}(\Omega))$},\\
    & \eps \nabla p_{\eps} \rightharpoonup \b0 \quad &&\text{weakly in $L^2(0,T;L^2(\Omega))$},\\
    & \eps \partial_t p_{\eps} \rightharpoonup 0 \quad &&\text{weakly in $L^2(0,T;W^{-1,2}(\Omega))$}.\\
\end{aligned}\end{equation*}
There exists a symmetric and traceless tensor function~$\bsigma$ such that
\begin{equation*}\begin{aligned}
    & \bsigma_{\eps} \rightharpoonup \bsigma \quad &&\text{weakly in $L^{\frac32}(0,T;L^{\frac32}(\Omega))$}.
\end{aligned}\end{equation*}
\end{proposition}

{\bf Proof -}
Weak convergence, up to a sub-sequence, for~$\u_{\eps}$, $\nabla\u_{\eps}$, $\partial_t \u_{\eps}$, $p_{\eps}$ and~$\bsigma_{\eps}$ follows from the bounds derived in Proposition~\ref{prop-bounds}.
Strong convergence on velocity~$\u_{\eps}$ can be deduced from the Aubin-Lions-Simon Lemma (see~\cite[p.102]{boyer12}).\\[0.2cm]
Since $p_\eps$ weakly converges to~$p$, we know that~$\partial_t p_\eps$ and~$\nabla p_\eps$ converge respectively to~$\partial_t p$ and~$\nabla p$ in the sense of distributions.
Thus, we deduce that~$\eps \partial_t p_\eps$ and~$\eps^{\frac{1}{2}}\nabla p_\eps$ converge to zero in the sense of distributions.
Since these two sequences are bounded in $L^2(0,T;W^{-1,2}(\Omega))$ and $L^2((0,T)\times\Omega)$, respectively, they converge to zero weakly in these spaces.\\[0.2cm]
The last point that remains to be proven is the non-negativity of the pressure~$p$.
For $\varphi \in L^3((0,T)\times \Omega)$, $\varphi\geq 0$, we have
\begin{equation}\label{eq:1342}
    \int_0^T \int_\Omega p_{\eps} \varphi
    =
    \underbrace{\iint_{p_{\eps}\geq 0} p_{\eps} \, \varphi}_{A}
    +
    \underbrace{\iint_{p_{\eps}<0} p_{\eps} \, \varphi}_{B}.
\end{equation}
First note that for all $\eps>0$ we have $A\geq 0$.
Moreover, using the Cauchy-Schwarz inequality and the energy estimate~\eqref{estk-1}, we obtain
\begin{equation*}
|B|
\leq \Big(\iint_{p_{\eps}<0} p_{\eps}^2 \Big)^{\frac12} \Big(\int_0^T \int_\Omega \varphi^2 \Big)^{\frac12}
\leq \eps^{\frac{1}{2}} \Big(\int_0^T \int_\Omega p_{\eps}V_\eps(p_\eps) \Big)^{\frac12} \Big(\int_0^T \int_\Omega \varphi^2 \Big)^{\frac12}
\xrightarrow{\eps\to 0} 0.
\end{equation*}
Finally, since $(p_{\eps})_{\eps>0}$ weakly converges to~$p$ in $L^{\frac32}((0,T)\times \Omega)$, passing to the limit as $\eps\to 0$, equality~\eqref{eq:1342} becomes
\begin{equation*}
\int_0^T \int_\Omega p \, \varphi \geq 0,
\end{equation*}
which corresponds to the non-negativity of the pressure~$p$.
\hfill$\blacksquare$

\vspace{0.3cm}

The aim at the end of this section is to demonstrate that we can pass to the limit in relations~\eqref{weak1k}, \eqref{weak2k} and~\eqref{weak3k}, especially by handling non-linear terms.

\paragraph{Step 1: Limit in~\eqref{weak1k}}
Given $\bvarphi\in L^3(0,T;W_0^{1,3}(\Omega))$, the convergences obtained in Proposition~\ref{prop:cv} indicate that
\begin{align}
& \langle \partial_t \u , \bvarphi \rangle - \int_\Omega p\, \div \, \bvarphi + \int_\Omega 2\overline{|\D\u|\D\u}:\D\bvarphi + \int_\Omega \bsigma:\S\bvarphi = \langle \f , \bvarphi \rangle, \label{weak1klim}
\end{align}
where the notation $\overline{|\D\u|\D\u}$ denotes the weak limit of $|\D\u_{\eps}||\D\u_{\eps}$ in $L^{\frac32}((0,T)\times \Omega)$.
\paragraph{Step 2: Limit in~\eqref{weak2k}}
In order to perform the limit $\eps\to 0$ in \eqref{weak2k}, we first note that
\[
\Big|\frac{|\S\u_{\eps}|^2}{|\S\u_{\eps}|+\eps} - |\S\u_{\eps}| \Big| =  \eps\frac{|\S\u_{\eps}|}{|\S\u_{\eps}|+\eps} \leq \eps \longrightarrow 0.
\]
Thus, for $\psi\in L^{\frac32}((0,T)\times \Omega)$ (in practice, this is done for more regular~$\psi_k$ functions, typically in $L^2(0,T;W_0^{1,2}(\Omega))$, and which converge to~$\psi$), when~$\eps$ goes to~$0$ in~\eqref{weak2k}, we get
\begin{align}
& \int_\Omega \psi \, \div \, \u - \int_\Omega 2\psi \overline{|\S\u|}  + \int_\Omega \psi \overline{V(p)} = 0, \label{weak2klim}
\end{align}
where $\overline{|\S\u|}$ and $\overline{V(p)}$ respectively correspond to the weak limits of~$|\S\u_{\eps}|$ and~$V_\eps(p_{\eps})$ in~$L^3((0,T)\times \Omega)$.
\paragraph{Step 3: Limit in~\eqref{weak3k}}
After integrating equation~\eqref{weak3k} with respect to time we obtain
\begin{align*}
&
\frac{1}{2}\|\u_{\eps}\|_2^2
+
\frac{\eps}{2}\|p_{\eps}\|_2^2
+
4 \int_0^T \|\D\u_{\eps}\|_3^3
+
\eps \int_0^T \|\nabla p_{\eps}\|_2^2
+
\int_0^T \int_\Omega p_{\eps}V_\eps(p_\eps)
= 
\frac{1}{2}\|\u_{\mathrm{init},\eps}\|_2^2
+
\int_0^T \langle \f_{\eps} , \u_{\eps} \rangle.
\end{align*}
In particular, we have
\begin{align}\label{eq:2214}
&
\frac{1}{2}\|\u_{\eps}\|_2^2
+
4 \int_0^T \|\D\u_{\eps}\|_3^3
+
\int_0^T \int_\Omega p_{\eps}V_\eps(p_\eps)
\leq
\frac{1}{2}\|\u_{\mathrm{init},\eps}\|_2^2
+
\int_0^T \langle \f_{\eps} , \u_{\eps} \rangle.
\end{align}
Since $(\u_\eps)_{\eps>0}$ converges weakly to~$\u$ in~$L^\infty(0,T;L^2(\Omega))$, we know that $\liminf_{\eps \to 0} \|\u_\eps\| \geq \|\u\|$.
We also use the strong convergence of the sequences~$(\u_{\mathrm{init},\eps})_{\eps>0}$ and~$(\f_{\eps})_{\eps>0}$ to perform the limit $\eps\to 0$ in~\eqref{eq:2214} and obtain
\begin{align}\label{weak3klim}
&
\frac{1}{2}\|\u\|_2^2
+
4 \int_0^T \int_\Omega \overline{|\D\u|^3}
+
\int_0^T \int_\Omega \overline{pV(p)}
\leq
\frac{1}{2}\|\u_{\mathrm{init}}\|_2^2
+
\int_0^T \langle \f , \u \rangle,
\end{align}
where $\overline{|\D\u|^3}$ and~$\overline{pV(p)}$ denote the weak limits of $(|\D\u_\eps|^3)_{\eps>0}$ and $(p_{\eps}V_\eps(p_\eps))_{\eps>0}$, respectively, in~$L^1((0,T)\times \Omega)$.

\subsection{Convergence of the non-linear terms and proof of the existence result}\label{subsec:limit-nonlinear}

In this subsection we will prove that
\begin{equation}\label{eq:1523}
\overline{|\D\u|\D\u} = |\D\u|\D\u,
\quad
\overline{V(p)} = \sqrt{p},
\quad
\overline{|\D\u|^3} = |\D\u|^3
\quad \text{and} \quad
\overline{|\S\u|} = |\S\u|.
\end{equation}
To this end, we make use of the following two lemmas:

\subsubsection{Convex analysis lemmas}\label{subsubsec:lemmas}

\begin{lemma}[Convexity argument]\label{lem:conv-faible}~\par
If $H:\mathcal D \subset \mathbb R^m \longrightarrow \mathbb R$ is a convex and lower semi-continuous function,\par
if $g_\eps: (0,T)\times \Omega \longrightarrow \mathcal D$ weakly converges to~$g$ in $L^q((0,T)\times \Omega)$ when~$\eps$ tends to~$0$,\par
if $H(g_\eps): (0,T)\times \Omega \longrightarrow \mathbb R$ weakly converges to~$\overline{H(g)}$ in $L^r((0,T)\times \Omega)$ when~$\eps$ tends to~$0$\par
then we have
$$H(g) \leq \overline{H(g)}.$$
\end{lemma}

{\bf Proof:}
This is a classical result that can be seen as a consequence of the equivalence of the notion of convex closed set for weak and strong topologies, see~\cite[p.\,38]{brezis83}.
For any $\varphi\in L^{r'}((0,T)\times \Omega)$, $\varphi \geq 0$ where $\frac{1}{r}+\frac{1}{r'}=1$, the function
$$\Phi:h\in L^q((0,T)\times \Omega) \longmapsto \int_0^T\int_\Omega H(h)\varphi$$
is convex and lower semi-continuous.
Following~\cite[Corollaire III.8]{brezis83} we deduce that if $g_\eps$ weakly converges to~$g$ in~$L^q((0,T)\times \Omega)^m$ then $\Phi(g) \leq \liminf_{\eps\to 0} \Phi(g_\eps)$.
We deduce that for all $\varphi\in L^{r'}((0,T)\times \Omega)$, $\varphi \geq 0$, we have
$$\int_0^T\int_\Omega H(g)\varphi \leq \int_0^T\int_\Omega \overline{H(g)}\varphi,$$
which corresponds to the announced result.
\hfill$\blacksquare$\\

Clearly, this result can also be extended to the case where the functions~$H$ is concave and upper semi-continuous. In this case, the function~$-H$ satisfies the assumptions of Lemma~\ref{lem:conv-faible} and we obtain
$$H(g) \geq \overline{H(g)}.$$

\begin{corollary}\label{cor:1}
We have the following inequalities
$$
\overline{pV(p)} \geq p\sqrt{p} \qquad \text{and} \qquad \overline{V(p)} \leq \sqrt{p}.
$$
\end{corollary}

{\bf Proof of the first inequality:}
Let $\varphi\in L^\infty((0,T)\times \Omega)$ with $\varphi\geq 0$. For all $\eps>0$ we have
    $$ \int_0^T\int_\Omega p_\eps V_\eps(p_\eps) \, \varphi 
    = \underbrace{\iint_{p_{\eps}\leq 0} p_\eps V_\eps(p_\eps) \, \varphi}_{A}
    + \underbrace{\iint_{p_{\eps}>0} p_\eps (V_\eps(p_\eps)-\sqrt{p_\eps}) \, \varphi}_{B}
    + \underbrace{\iint_{p_{\eps}>0} p_\eps \sqrt{p_\eps} \, \varphi}_{C}.$$
    \begin{itemize}
        \item[--] Since for all $x\leq 0$ we have $xV_\eps(x)\geq 0$ we deduce that $A\geq 0$.
        \item[--] Since for all $x>0$ we have $|V_\eps(x)-\sqrt{x}|\leq \eps$, we deduce that $\lim_{\eps\to 0} B=0$.
        \item[--] Finally, applying Lemma~\ref{lem:conv-faible} with the convex and continuous function $H:x\in \mathbb R_+ \longmapsto x\sqrt{x}$, and using the fact that the limit pressure~$p$ is non-negative, we deduce that $\liminf_{\eps\to 0} C \geq \int_0^T\int_\Omega p \sqrt{p} \, \varphi.$
    \end{itemize}
    We then conclude that 
    $$\int_0^T\int_\Omega \overline{p V(p)} \, \varphi \geq \int_0^T\int_\Omega p \sqrt{p} \, \varphi.$$
    
{\bf Proof of the second inequality:}
Let $\varphi\in L^{\frac{3}{2}}((0,T)\times \Omega)$ with $\varphi\geq 0$. For all $\eps>0$ we have
    $$\int_0^T\int_\Omega V_\eps(p_\eps) \, \varphi 
    = \underbrace{\iint_{p_{\eps}\leq 0} V_\eps(p_\eps) \, \varphi}_{A}
    + \underbrace{\iint_{p_{\eps}>0} V_\eps(p_\eps) \, \varphi}_{B}.$$
    \begin{itemize}
        \item[--] Since for all $x\leq 0$ we have $V_\eps(x)\leq 0$, we deduce $A\leq 0$.
        \item[--] Since for all $x> 0$ we have $V_\eps(x)\leq \sqrt{x}$, we deduce $B\leq \iint_{p_{\eps}>0} \sqrt{p_\eps} \, \varphi$. Then, applying Lemma~\ref{lem:conv-faible} (in its concave version, as detailed immediately following the lemma's proof) with the concave and continuous function $H:x\in \mathbb R_+ \longmapsto \sqrt{x}$, and using the fact that the limit pressure~$p$ is non-negative, we deduce that $\liminf_{\eps\to 0} B \leq \int_0^T\int_\Omega \sqrt{p} \, \varphi.$
    \end{itemize}
    We then conclude that 
    $$\int_0^T\int_\Omega \overline{V(p)} \, \varphi \leq \int_0^T\int_\Omega \sqrt{p} \, \varphi,$$
    that completes the proof of Corollary~\ref{cor:1}.
    \hfill$\blacksquare$\\

The following lemma can be found in~\cite[Lemma 1]{zhikov10}, see also~\cite[p.1715 and p.1730]{fang22}:
\begin{lemma}\label{lemFangGuo22}
If $\D\u_{\eps} \rightharpoonup \D\u$ weakly in $L^3((0,T)\times \Omega)$ and $|\D\u_{\eps}|\D\u_{\eps} \rightharpoonup \overline{|\D\u|\D\u}$ weakly in $L^\frac32((0,T)\times \Omega)$ then
\[
\int_0^T\int_\Omega \Big( \overline{|\D\u|\D\u}:\D\u - 2\overline{|\D\u|^3} \Big) \leq 0.
\]
where $\overline{|\D\u|^3}$ refers to the weak limit of $|\D\u_\eps|^3$ in $L^1((0,T)\times \Omega)$.\\
Furthermore, if there is equality, then $\overline{|\D\u|\D\u} = |\D\u|\D\u$.
\end{lemma}

\subsubsection{Existence proof}\label{subsubsec:existence}

With the help of Lemmas~\ref{lem:conv-faible} and~\ref{lemFangGuo22}, we can conclude the existence proof stated in the theorem by choosing $\bvarphi=\u$ in equation~\eqref{weak1klim} and $\psi=p$ in equation~\eqref{weak2klim}, integrating both equations over time, and summing the resulting expressions.
This yields
\begin{align*}
\frac{1}{2}\|\u\|_2^2
+
\int_0^T\int_\Omega 2\overline{|\D\u|\D\u}:\D\u
+
\int_0^T\int_\Omega p \overline{V(p)}
+
\int_0^T \int_\Omega \Big( \bsigma:\S\u - 2p\overline{|\S\u|} \Big) 
=
\frac{1}{2}\|\u_{\mathrm{init}}\|_2^2
+
\int_0^T \langle \f , \u \rangle.
\end{align*}
By subtracting equation~\eqref{weak3klim}, we obtain
\begin{align}\label{eq1134}
\underbrace{\int_0^T\int_\Omega 2\Big( \overline{|\D\u|\D\u}:\D\u - 2\overline{|\D\u|^3} \Big)}_{A}
+
\underbrace{\int_0^T\int_\Omega \Big( p\overline{V(p)} - \overline{pV(p)} \Big)}_{B}
+
\underbrace{\int_0^T \int_\Omega \Big( \bsigma:\S\u - 2p\overline{|\S\u|} \Big)}_{C}
\geq 0.
\end{align}
\begin{itemize}
\item[--] By Lemma~\ref{lemFangGuo22}, we have $A\leq 0$.
\item[--] From Corollary~\ref{cor:1}, we directly obtain $B\leq 0$.
\item[--] By the convexity of $x\mapsto |x|$ and Lemma~\ref{lem:conv-faible}, we have $|\S\u| \leq \overline{|\S\u|}$. Moreover, we have $|\bsigma|\leq p$. Consequently, by applying the Cauchy-Schwarz inequality we write
\begin{equation}\label{eq:1208}
    \bsigma:\S\u \leq 2|\bsigma||\S\u| \leq 2p|\S\u| \leq 2p\overline{|\S\u|}.
\end{equation}
which implies that $C\leq 0$.
\end{itemize}
Since $A+B+C\geq 0$ with $A\leq 0$, $B\leq 0$ and $C\leq 0$, it follows that $A=B=C=0$.
\begin{itemize}
\item[--] According to the conclusion of Lemma~\ref{lemFangGuo22}, the equality~$A=0$ implies $\overline{|\D\u|\D\u} = |\D\u|\D\u$.
\item[--] Similarly, the equality $B=0$ implies $\overline{V(p)}=\sqrt{p}$.
\item[--] Finally, even if this point is not useful in the proof of Theorem~\ref{th:1}, we can note that~$C=0$ implies $\bsigma:\S\u = 2p\overline{|\S\u|}$ so that equality holds  in~\eqref{eq:1208}, implying that $\bsigma:\S\u = 2p|\S\u|$.
\end{itemize}
In particular, this last point proves that $2p|\S\u|=2p\overline{|\S\u|}$ but does not prove that $|\S\u|=\overline{|\S\u|}$.
To obtain the latter result, we proceed in the following four steps.\\[0.2cm]
{\bf Step 1:} $\overline{|\D\u|^3} = |\D\u|^3$\\
This is a direct consequence of Lemma~\ref{lemFangGuo22}. Indeed, since $A=0$ and $\overline{|\D\u|\D\u} = |\D\u|\D\u$, we have
\[
\int_0^T\int_\Omega \Big( |\D\u|^3 - \overline{|\D\u|^3} \Big) = 0.
\]
Moreover the function $x\mapsto x^3$ being convex, we apply Lemma~\ref{lem:conv-faible} to obtain $|\D\u|^3 - \overline{|\D\u|^3} \leq 0$.
We deduce the equality $\overline{|\D\u|^3} = |\D\u|^3$.\\[0.2cm]
{\bf Step 2:} $\overline{|\D\u|^2} = |\D\u|^2$\\
Using the convexity of $x\mapsto x^2$ and the concavity of $x\mapsto x^{\frac{2}{3}}$, we deduce from Lemma~\ref{lem:conv-faible} that
\begin{equation*}
    |\D\u|^2 
    \leq \overline{|\D\u|^2}
    = \overline{(|\D\u|^3)^{\frac{2}{3}}}
    \leq \overline{|\D\u|^3}^{\frac{2}{3}}
    = |\D\u|^2.
\end{equation*}
All these terms are therefore equal, and in particular $\overline{|\D\u|^2} = |\D\u|^2$.\\[0.2cm]
{\bf Step 3:} $\overline{|\S\u|^2} = |\S\u|^2$\\
We use the identity
\begin{equation}\label{eq:1317}
|\D\u_\eps|^2 = |\S\u_\eps|^2 + \frac{1}{6}|\div\, \u_\eps|^2.
\end{equation}
Let $\xi\geq 0$ and $\zeta\geq 0$ be the elements of $L^\frac{3}{2}((0,T)\times \Omega)$ such that $\overline{|\S\u|^2} = |\S\u|^2 + \xi$ and $\overline{|\div\, \u|^2} = |\div\, \u|^2 + \zeta$ then
we pass to the limit in~\eqref{eq:1317} using the fact that $\overline{|\D\u|^2} = |\D\u|^2$:
\begin{equation*}
|\D\u|^2 = |\S\u|^2 + \xi + \frac{1}{6}|\div\, \u|^2 + \frac{1}{6}\zeta.
\end{equation*}
Since we also have $|\D\u|^2 = |\S\u|^2 + \frac{1}{6}|\div\, \u|^2$, it follows that $\xi + \frac{1}{6}\zeta=0$. Given the positivity of $\xi$ and $\zeta$, we conclude that $\xi=\zeta=0$, and in particular $\overline{|\S\u|^2} = |\S\u|^2$.\\[0.2cm]
{\bf Step 4:} $\overline{|\S\u|} = |\S\u|$\\
Using the convexity of $x\mapsto |x|$, and the concavity of $x\mapsto x^{\frac{1}{2}}$, we deduce from Lemma~\ref{lem:conv-faible} that
\begin{equation*}
    |\S\u|
    \leq \overline{|\S\u|}
    = \overline{(|\S\u|^2)^{\frac{1}{2}}}
    \leq \overline{|\S\u|^2}^{\frac{1}{2}}
    = |\S\u|.
\end{equation*}
All these terms are therefore equal, and in particular $\overline{|\S\u|} = |\S\u|$.\\[0.2cm]
To conclude the proof of Theorem~\ref{th:1}, simply note that equations~\eqref{weak1klim}, \eqref{weak2klim} and~\eqref{weak3klim} exactly correspond to Definition~\ref{definition} of a weak solution since we have proved~\eqref{eq:1523}.
\hfill$\blacksquare$

\subsection{Proof of the uniqueness result}\label{subsec:uniqueness}

The proof of the uniqueness of a weak solution is classical. Consider two weak solutions $(\u_1,p_1,\bsigma_1)$ and $(\u_2,p_2,\bsigma_2)$. We write the equations~\eqref{weak1} and~\eqref{weak2} for both solutions, and then take the difference. Denoting $\u=\u_2-\u_1$, $p=p_2-p_1$ and $\bsigma=\bsigma_2-\bsigma_1$, we obtain, for all $\bvarphi\in L^3(0,T;W_0^{1,3}(\Omega))$ and $\psi\in L^{\frac32}((0,T)\times \Omega)$:
\begin{align*}
& \langle \partial_t \u , \bvarphi \rangle - \int_\Omega p\, \div \, \bvarphi + \int_\Omega 2(|\D\u_2|\D\u_2-|\D\u_1|\D\u_1):\D\bvarphi + \int_\Omega \bsigma:\S\bvarphi = 0, \\
& \int_\Omega \psi \, \div \, \u - \int_\Omega 2\psi (|\S\u_2|-|\S\u_1|) + \int_\Omega \psi (\sqrt{p_2}-\sqrt{p_1}) = 0.
\end{align*}
We then choose $\bvarphi=\u$ and $\psi=p$, and add the resulting equations:
\begin{equation}\label{eq:2043}
\begin{aligned}
    \frac{1}{2}\frac{d}{dt}\int_\Omega |\u|^2 
    & + \int_\Omega \underbrace{(|\D\u_2|\D\u_2-|\D\u_1|\D\u_1):(\D\u_2-\D\u_1)}_{A} \\
&     + \int_\Omega \underbrace{(\sqrt{p_2}-\sqrt{p_1})(p_2-p_1)}_{B} \\
&     + \int_\Omega \underbrace{(\bsigma_2-\bsigma_1):(\S\u_2-\S\u_1) - 2(p_2-p_1)(|\S\u_2|-|\S\u_1|)}_{C} \leq 0.
\end{aligned}
\end{equation}
We remark that $A\geq (|\D\u_2|+|\D\u_1|)(|\D\u_2|-|\D\u_1|)^2 \geq 0$ and that $B= (\sqrt{p_2}+\sqrt{p_1})(\sqrt{p_2}-\sqrt{p_1})^2 \geq 0$. Moreover, using Proposition~\ref{prop:1} and the Cauchy-Schwarz inequality, we also have $C\geq 0$.
Consequently, estimate~\eqref{eq:2043} implies that $\|\u\|_2$ is a positive and non-increasing function of time.
Since~$\u_2$ and~$\u_1$ satisfy the same initial condition, we deduce that $\u=\b0$.\\
We next deduce from the estimate~\eqref{eq:2043} that $A=B=C=0$.
In particular, the equality~$B=0$ implies that $p=0$.
\hfill$\blacksquare$

\section{Adaptation to more general models}\label{sec:generalisation}

\subsection{Taking into account the variation in volume fraction}\label{subsec:phi-variation}

As explained in the introduction, the complete granular flow model also describes the evolution of the grain volume fraction, denoted~$\phi$.
Although the previous results have overcome a number of difficulties relating to complex flows (non-linear rheology and dilatation phenomena), the complete model given by~\eqref{syst1}--\eqref{syst0} in the introduction presents additional difficulties.
The key point is to specify the functions~$\alpha$ and~$\beta$ which depend on the volume fraction~$\phi$. More specifically, it is important to understand how to describe the expansion law expressing $\div \, \u$.
From Roux-Radjaï's work~\cite{roux98}, we know that this law can be written as
$$\div \, \u = 2K |\S\u|(I-I_{\mathrm{eq}}(\phi))
\qquad \text{where} \quad I=\frac{2d|\S\u|}{\sqrt{p/\rho_0}}.$$
The non-dimensional number~$I$ denotes the inertial number (depending on the diameter~$d$ of the grains, and on their density~$\rho_0$) and $I_{\mathrm{eq}}(\phi)$ its equilibrium value. For the latter, several examples are given in the literature, see~\cite{breard22, jop_etal,gdrmidi,robinson23}, and we choose~\cite[p.929]{schaeffer19}:
$$I_{\mathrm{eq}}(\phi) = \frac{\phi_{\mathrm{max}}-\phi}{\phi-\phi_{\mathrm{min}}},$$
where the constants $0<\phi_{\mathrm{min}} < \phi_{\mathrm{max}}<1$ corresponds to the "extreme" values of the volume fraction~$\phi$.
This choice will ensure mathematically that $\phi$ remains effectively bounded between~$\phi_{\mathrm{min}}$ and~$\phi_{\mathrm{max}}$, see Proposition~\ref{prop:bound-infini-phi} below.
In order to ensure stability of the model (see~\cite{barker23,barker17,barker15,Chupin-Dubois24}), we choose the proportionality coefficient as $K=\frac{\phi-\phi_{\mathrm{min}}}{\delta \phi \, I}$ with $\delta \phi=\phi_{\mathrm{max}}-\phi_{\mathrm{min}}$, so that the dilatation law is written as
$$\div \, \u =  2 \frac{\phi-\phi_{\mathrm{min}}}{\delta \phi} |\S\u| - \frac{\phi_{\mathrm{max}}-\phi}{\delta \phi} \frac{\sqrt{p}}{d\sqrt{\rho_0}}.$$
Regarding the mathematical analysis, we propose here a first theoretical result on the existence and uniqueness for a model that is relatively close to the full model.
This considered model is as follows:
\begin{numcases}{}
\label{systphi1}
\mathcal D_t(\phi,\u) + \nabla p - \div ( 2|\D\u|\D\u) = \f + \div \, \bsigma, \\
\label{systphi3}
\bsigma:\S\u=2 (\phi-\phi_{\mathrm{min}}) p|\S\u|, \quad |\bsigma| \leq (\phi-\phi_{\mathrm{min}}) p, \quad \bsigma=\bsigma^{\intercal} \quad\text{and}\quad \mathrm{tr}\, \bsigma=0, \\
\label{systphi4}
\div \, \u = 2 (\phi-\phi_{\mathrm{min}}) |\S\u| - (\phi_{\mathrm{max}}-\phi) \sqrt{p}, \\
\label{systphi0}
\partial_t \phi + \div(\phi\u) = \xi \mathcal H,
\end{numcases}
where $$\mathcal D_t(\phi,\u) = \frac{1}{2}\Big( \big(
\partial_t (\phi\u) + \div(\phi\u\otimes \u) \big)
+
\big(
\phi(\partial_t\u+\u\cdot \nabla \u)
\big)\Big).$$
In practice, the only difference with the full model lies in the additional term $\xi \mathcal H$ where~$\xi$ is a real parameter, fixed and positive (possibly very small).
This term~$\mathcal{H}$ will be chosen in the following form
\begin{equation*}
    \mathcal{H} = \Delta \phi - \phi\sqrt{p}.
\end{equation*}
The contribution of $\Delta \phi$ ensures compactness for sequence~$(\phi_\eps)_{\eps>0}$ satisfying such a problem, and gives sense to the non-linear terms such as $(\phi_{\mathrm{max}}-\phi) \sqrt{p}$.
Furthermore, the contribution of~$\phi\sqrt{p}$ ensures that~$\phi\leq \phi_{\mathrm{max}}-\xi$, thereby allowing the energy estimate to yield a bound on the pressure.
\begin{remark}
The form of the time derivative as the average between the convected derivative and its conservative form ensures an energy estimate independent of the evolution equation~\eqref{systphi0} on~$\phi$.
In the case where $\xi=0$, both forms are equivalent.
At this stage, we do not know how to obtain sufficient regularity for the solutions to assert the existence of a solution when~$\xi=0$.
In particular, we do not know whether the solutions constructed here converge to a solution when~$\xi$ goes to~$0$.
\end{remark}

\begin{theorem}\label{th:2}
Let $\xi>0$, $\u_{\mathrm{init}}\in L^2(\Omega)$, $\phi_{\mathrm{init}}\in L^2(\Omega)$ and $\f\in L^{\frac32}(0,T;W^{-1,\frac32}(\Omega))$.\\
If $\phi_{\mathrm{min}}\leq \phi_{\mathrm{init}} \leq \phi_{\mathrm{max}}-\xi$ then there exists a weak solution of~\eqref{systphi1}--\eqref{systphi0} satisfying the initial conditions $\u|_{t=0}=\u_{\mathrm{init}}$ and~$\phi|_{t=0}=\phi_{\mathrm{init}}$.
\end{theorem}

The proof follows the same arguments as in the case of Theorem~\ref{th:1}. The only difference lies in the control of~$\phi$.
We must therefore ensure that there is always an energy estimate that, in addition to providing control over~$\u$ and~$p$, also allows us to control~$\phi$. 
To avoid repeating the full proof of the previous theorem, we simply derive the energy estimate in a formal manner.

\begin{proposition}\label{prop:2}
Any regular solution $(\u,p,\bsigma,\phi)$ to~\eqref{systphi1}--\eqref{systphi0} satisfies
\begin{equation}\label{eq:2115}
    \frac{d}{dt}\int_\Omega \phi\frac{|\u|^2}{2} + 4 \int_\Omega |\D\u|^3 + \int_\Omega (\phi_{\mathrm{max}}-\phi) p\sqrt{p} = \int_\Omega \f\cdot \u.
\end{equation}
\end{proposition}
{\bf Proof -}
The result essentially follows by taking the scalar product of equation~\eqref{systphi1} with~$\u$.
We note that, independently of equation~\eqref{systphi0}, we have
$$\mathcal D_t(\phi,\u)\cdot \u = \partial_t \Big( \phi\frac{|\u|^2}{2} \Big) + \div \Big( \phi\frac{|\u|^2}{2}\u \Big).
$$
Equation~\eqref{systphi4} allows us to control the pressure term
$$
\nabla p\cdot \u = \div(p\u) - p\div\,\u = \div(p\u) -2(\phi-\phi_{\mathrm{min}})p|\S\u| + (\phi_{\mathrm{max}}-\phi) p\sqrt{p},
$$
while equation~\eqref{systphi3} is useful for cancelling out the stress part:
$$
\div \, \bsigma \cdot \u = \div(\bsigma \cdot \u) - \bsigma : \S\u = \div(\bsigma \cdot \u) - 2(\phi-\phi_{\mathrm{min}})p|\S\u|.
$$
After integration over $\Omega$, the sum of these contributions provides the result stated in Proposition~\ref{prop:2}.
\hfill$\blacksquare$\\[0.2cm]
To use equation~\eqref{eq:2115}, we need to prove that the terms containing $\phi|\u|^2$ and~$(\phi_{\mathrm{max}}-\phi) p\sqrt{p}$ are positive, more precisely that $0< \phi \leq \phi_{\mathrm{max}}$.
In addition, if we want to obtain an estimate on the pressure, we need to ensure that $\phi_{\mathrm{max}}-\phi$ does not vanish.
\begin{proposition}\label{prop:bound-infini-phi}
Any regular solution $(\u,p,\bsigma,\phi)$ to~\eqref{systphi1}--\eqref{systphi0} such that $\phi_{\mathrm{min}}\leq \phi|_{t=0} \leq \phi_{\mathrm{max}}-\xi$ satisfies 
\begin{equation}\label{estimationLinfty}
\phi_{\mathrm{min}} \leq \phi \leq \phi_{\mathrm{max}}-\xi,
\end{equation}
and
\begin{equation}\label{eq:2116}
\frac{d}{dt}\|\phi\|_2^2 + 2\xi \|\nabla \phi\|_2^2 \leq \phi_{\mathrm{max}}^2 \int_\Omega (|\div\, \u| + 2\sqrt{p}).
\end{equation}
\end{proposition}

{\bf Proof -- Step 1: Upper bound}
First note that equation~\eqref{systphi0} implies that, for $\beta:\mathbb R\longrightarrow\mathbb R$, we have
\begin{equation}\label{dt(beta)}
    \partial_t \beta(\phi) + \div ( \beta(\phi) \u ) + (\phi \beta'(\phi) - \beta(\phi)) \div \, \u = \xi \beta'(\phi)\Delta \phi - \xi\phi\beta'(\phi)\sqrt{p},
\end{equation}
In order to achieve the upper bound for~$\phi$, we choose the function~$\beta$ defined by
\begin{equation*}
\beta(\phi)=\left\{\begin{aligned}
& \hspace{1.3cm} 0 \quad && \text{if $\phi<\phi_{\mathrm{max}}-\xi$},\\
& \phi-(\phi_{\mathrm{max}}-\xi) \quad && \text{if $\phi \geq \phi_{\mathrm{max}}-\xi$}.
\end{aligned}\right.
\end{equation*}
Integrating~\eqref{dt(beta)} over~$\Omega$, and thanks to this choice of the function~$\beta$, we obtain
\begin{equation}\label{dt(int(beta))}
    \frac{d}{dt} \int_\Omega \beta(\phi) + \int_{\mathcal E_+} (\phi_{\mathrm{max}}-\xi) \div \, \u = \xi \int_{\mathcal E_+} \Delta \phi - \int_{\mathcal E_+} \xi \phi\sqrt{p},
\end{equation}
where $\mathcal E_+ = \{x\in \Omega~;~ \phi(t,x) \geq \phi_{\mathrm{max}}-\xi\}$.
The unit outgoing normal vector at~$\mathcal E_+$ is given by $-\nabla \phi/|\nabla \phi|$ so that Stokes' formula allows us to write
$$\xi \int_{\mathcal E_+} \Delta \phi
= - \xi \int_{\partial \mathcal E_+} |\nabla \phi| \leq 0.$$
Moreover, the expression for the divergence of the velocity given by equation~\eqref{systphi4} indicates that
$$\phi \geq \phi_{\mathrm{max}}-\xi \quad \Longrightarrow \quad \div \, \u + \xi\sqrt{p} \geq 0$$
We then deduce that
$$\int_{\mathcal E_+} (\phi_{\mathrm{max}}-\xi) \div \, \u + \int_{\mathcal E_+} \xi \phi\sqrt{p}
\geq \int_{\mathcal E_+} (\phi_{\mathrm{max}}-\xi) (\div \, \u + \xi \sqrt{p})
\geq 0.$$
Consequently, equation~\eqref{dt(int(beta))} implies $\frac{d}{dt} \int \beta(\phi) \leq 0$.
Thus, if $\phi|_{t=0} \leq \phi_{\mathrm{max}}-\xi$, i.e. $\beta(\phi|_{t=0})=0$, then $\beta(\phi)=0$ so the result follows, namely $\phi \leq \phi_{\mathrm{max}}-\xi$.\\[0.2cm]
{\bf Step 2: Lower bound}
In the same spirit, to obtain the lower bound, we use equation~\eqref{dt(beta)} together with the function~$\beta$, defined by 
\begin{equation*}
\beta(\phi)=\left\{\begin{aligned}
& \hspace{0.6cm} 0 \quad && \text{if $\phi > \phi_{\mathrm{min}}$},\\
& \phi_{\mathrm{min}}-\phi \quad && \text{if $\phi \leq \phi_{\mathrm{min}}$}.
\end{aligned}\right.
\end{equation*}
By integrating \eqref{dt(beta)} over~$\Omega$, we get
\begin{equation}\label{dt(int(beta))1}
    \frac{d}{dt} \int_\Omega \beta(\phi) - \int_{\mathcal E_-} \phi_{\mathrm{min}} \div \, \u = -\xi \int_{\mathcal E_-} \Delta \phi + \int_{\mathcal E_-} \xi \phi\sqrt{p},
\end{equation}
where $\mathcal E_- = \{x\in \Omega~;~ \phi(t,x)\leq\phi_{\mathrm{min}}\}$.
This time, the unit outgoing normal vector at~$\mathcal E_-$ is given by $\nabla \phi/|\nabla \phi|$ so that Stokes' formula allows us to write 
$$- \xi \int_{\mathcal E_-} \Delta \phi
= - \xi \int_{\partial \mathcal E_-} |\nabla \phi| \leq 0.$$
The expression for the divergence of the velocity in equation~\eqref{systphi4} gives
$$
\phi \leq \phi_{\mathrm{min}} \quad \Longrightarrow \quad \div \, \u \leq - (\phi_{\mathrm{max}}-\phi_{\mathrm{min}}) \sqrt{p}.
$$
Since $\xi\leq \phi_{\mathrm{max}}-\phi_{\mathrm{min}}$, we obtain the following inequality
$$
\int_{\mathcal E_-} \xi \phi\sqrt{p}
+
\int_{\mathcal E_-} \phi_{\mathrm{min}} \div \, \u
\leq \int_{\mathcal E_-} \phi_{\mathrm{min}} (\xi - (\phi_{\mathrm{max}}-\phi_{\mathrm{min}}) )\sqrt{p}
\leq 0.
$$
Consequently, equation~\eqref{dt(int(beta))1} implies that $\frac{d}{dt} \int \beta(\phi) \leq 0$.
Thus, if we have $\phi|_{t=0} \geq \phi_{\mathrm{min}}$, i.e. $\beta(\phi|_{t=0})=0$, then $\beta(\phi)=0$ so that the result follows, namely $\phi \geq \phi_{\mathrm{min}}$.\\[0.2cm]
{\bf Step 3: $H^1$ estimate}
By multiplying equation~\eqref{systphi0} by~$2\phi$ and integrating over~$\Omega$, we obtain
\begin{equation*}
    \frac{d}{dt}\int_\Omega |\phi|^2 + 2\xi \int_\Omega |\nabla \phi|^2 = -\int_\Omega 2\phi \, \div(\phi\u) - \int_\Omega 2\phi^2\sqrt{p}.
\end{equation*}
Using integration by parts, we can write the right-hand side in the form
\begin{equation*}
    \frac{d}{dt}\|\phi\|_2^2 + 2\xi \|\nabla \phi\|_2^2 = -\int_\Omega \phi^2 (\div\, \u + 2\sqrt{p}).
\end{equation*}
The bound $|\phi|\leq \phi_{\mathrm{max}}$ allows to conclude.
\hfill$\blacksquare$

\paragraph{Ideas for the proof of Theorem~\ref{th:2}}

The method is the same as the one described in the proof of Theorem~\ref{th:1}. We construct a sequence of solutions $(\u_\eps,p_\eps,\bsigma_\eps,\phi_\eps)_{\eps>0}$ of an approximate problem, and we show that this solution converges to the solution of the problem~\eqref{systphi1}--\eqref{systphi0} when~$\eps$ tends to~$0$.\\[0.2cm]
We focus on estimating the volume fraction~$\phi_\eps$, since the other quantities are bounded similarly to those in the proof of Theorem~\ref{th:1}. The key point is to ensure that the sequence~$(\phi_\eps)_{\eps>0}$ converges strongly to a solution~$\phi$, which will enable us to pass to the limit in all terms.\\[0.2cm]
Given that $\phi_{\mathrm{max}} - \phi_\eps \geq \xi > 0$, estimate~\eqref{eq:2115} ensures that the sequences $(\div\, \u_\eps)_{\eps>0}$ and~$(\sqrt{p_\eps})_{\eps>0}$ are uniformly bounded in~$L^3((0,T)\times \Omega)$ with respect to~$\eps$.
From estimate~\eqref{eq:2116}, we deduce that $(\phi_\eps)_{\eps>0}$ is uniformly bounded in $L^2(0,T;H^1(\Omega))$.\\[0.2cm]
To obtain compactness, we consider the sequence $(\partial_t\phi_\eps)_{\eps>0}$.
By writing
\begin{equation*}
    \partial_t \phi = \xi \Delta \phi - \div(\phi\u) - \phi\sqrt{p},
\end{equation*}
we deduce from the previous estimates that $(\partial_t\phi_\eps)_{\eps>0}$ is uniformly bounded in~$L^2(0,T;H^{-1}(\Omega))$. Finally, the Aubin-Lions-Simon Lemma (see~\cite[p.102]{boyer12}) leads to the conclusion that $(\phi_\eps)_{\eps>0}$ converges strongly to~$\phi$ in $L^2((0,T)\times \Omega)$.
\hfill$\blacksquare$

\begin{remark}
Propositions~\ref{prop:2} and~\ref{prop:bound-infini-phi} are true even if $\xi=0$.
However, the proof of Theorem~\ref{th:2} is no longer correct.
Indeed, when $\xi=0$, we don't have $H^1$ estimate for~$\phi$ (see relation~\eqref{eq:2116}), nor an estimate for~$p$ (see relation~\eqref{eq:2115}, especially when $\phi_{\mathrm{max}}-\phi$ vanishes).
To prove such a theorem when $\xi=0$, we therefore need to find other arguments.
\end{remark}

\subsection{\texorpdfstring{$\mu(I)$}{mu\texttwosuperior}-rheology}\label{subsec:mu(I)}

The $\mu(I)$-rheology relates stress, pressure and shear in a manner similar to equation~\eqref{syst-loi-seuil}.
The main difference is that the stress threshold, corresponding to~$\alpha(\phi)p$ in~\eqref{syst-loi-seuil}, is written~$\mu(I)p$, the function~$\mu$ being experimentally given with respect to the inertial number~$I$, see for example~\cite{andreotti12}.
In this context, if Barker's stability conditions~\cite{schaeffer19} are to be satisfied, the dilation relation is also imposed, and involves $I$-dependent functions.
By following these principles the following model is proposed in~\cite{Chupin-Dubois24}:
\begin{numcases}{}
\label{syst1bis} 
\phi \rho_0 \big( \partial_t \u + \u\cdot \nabla \u \big) + \nabla p = \phi \rho_0 \g + \div \, \bsigma, \\
\label{syst3bis}
\bsigma:\S\u=2\mu(I)p|\S\u|, \quad |\bsigma| \leq \mu(I) p, \quad \bsigma=\bsigma^{\intercal} \quad\text{and}\quad \mathrm{tr}\, \bsigma=0, \\
\label{syst4is}
\div \, \u = 2 F(I) |\S\u| - \gamma I_{\mathrm{eq}}(\phi)F(I_{\mathrm{eq}}(\phi))\sqrt{p}, \\
\label{syst0bis}
\partial_t \phi + \div(\phi\u) = 0,
\end{numcases}
where $\displaystyle \gamma=2/d\sqrt{\rho_0}$ and where the dimensionless functions $\mu$ and $F$ are explicitly given (see~\cite{Chupin-Dubois24} for more details). In particular, these functions are smooth, continuous, and equal at~$0$. We denote $\alpha_0=\mu(0)=F(0)$.

\paragraph{Dimensionless procedure}

In order to rewrite the system in dimensionless form, we apply the following change of variables and unknowns:
\[t=T\widetilde t, \quad x=L\widetilde x, \quad \u=U\widetilde{\u}, \quad p=gL\rho_0\widetilde p, \quad \bsigma = gL\rho_0\widetilde{\bsigma},\]
where $T$, $L$ and $U$ represent the characteristic time, length and velocity scales, respectively.
Introducing the dimensionless numbers~$\eps$, $\mathfrak{Fr}$ and~$\mathfrak{Di}$ defined as follows
\[\eps=\frac{TU}{L}, \quad \mathfrak{Fr}^2=\frac{U^2}{gL}, \quad \mathfrak{Di}=\frac{d}{L}\mathfrak{Fr},\]
yields the following system where the tilde notation has been omitted for conciseness:
\begin{numcases}{}
\label{syst1ter} 
\phi \big( \partial_t \u + \eps \u\cdot \nabla \u \big) + \frac{\eps}{\mathfrak{Fr}^2}\nabla p = -\frac{\eps}{\mathfrak{Fr}^2} \phi \e + \frac{\eps}{\mathfrak{Fr}^2}\div \, \bsigma, \\
\label{syst3ter}
\bsigma:\S\u=2\mu(\mathfrak{Di}\, I)p|\S\u|, \quad |\bsigma| \leq \mu(\mathfrak{Di}\, I) p, \quad \bsigma=\bsigma^{\intercal} \quad\text{and}\quad \mathrm{tr}\, \bsigma=0, \\
\label{syst4ter}
\div \, \u = 2 F(\mathfrak{Di}\, I) |\S\u| - \frac{1}{\mathfrak{Di}} I_{\mathrm{eq}}(\phi)F(I_{\mathrm{eq}}(\phi)) \sqrt{p}, \\
\label{syst0ter}
\partial_t \phi + \eps \div(\phi\u) = 0.
\end{numcases}
Note that in this dimensionless form, the inertial number~$I$ is written as:~$\displaystyle I=\frac{2|\S\u|}{\sqrt{p}}$.

\paragraph{Asymptotic problem}

We aim to study this system when the volume fraction~$\phi$ is nearly constant and close to the maximum packing fraction $\phi_{\mathrm{max}}$.
Choosing~$\eps$ as a small parameter, let us introduce
$$\phi = \phi_0+\eps \psi,$$
where~$\phi_0$ is a constant close to~$\phi_{\mathrm{max}}$ in the following sense: $\phi_{\mathrm{max}}-\phi_0 = \mathcal O(\eps)$.
More generally, we write
\[\u=\v+\mathcal O(\eps), \quad p=q+\mathcal O(\eps) \quad \text{and} \quad \bsigma=\btau + \mathcal O(\eps).\]
To retain both rheological and expansion-related contributions, we also assume that
\[\mathfrak{Fr}^2 = \mathcal O(\eps) \quad \text{and} \quad \mathfrak{Di} = \mathcal O(\eps).\]
A typical illustration of this situation is shown below:
$$
L=10^{-1}\, \mathrm{m}, \quad U=10^{-1}\, \mathrm{m.s}^{-1},\quad T=10^{-2}\, \mathrm{s} \quad d=10^{-2}\, \mathrm{m} \quad \text{and} \quad g=10\, \mathrm{m.s}^{-2}.
$$
In this case, we have $\eps=10^{-2}$, $\mathfrak{Fr}^2=10^{-2}$ and~$\mathfrak{Di}=10^{-2}$.
By keeping only the terms of leading order in~$\eps$ in the system~\eqref{syst1ter}--\eqref{syst0ter}, we obtain
\begin{numcases}{}
\label{syst1qua} 
\phi_0 \partial_t \v + \lambda_0\nabla q = -\phi_0 \lambda_0 \e + \lambda_0 \div \, \btau, \\
\label{syst3qua}
\btau:\S\v=2\alpha_0 q|\S\v|, \quad |\btau| \leq \alpha_0 q, \quad \btau=\btau^{\intercal} \quad\text{and}\quad \mathrm{tr}\, \btau=0, \\
\label{syst4qua}
\div \, \v = 2 \alpha_0 |\S\v| - \gamma_0 \sqrt{q}, \\
\label{syst0qua}
\partial_t \psi + \phi_0 \div\, \v = 0,
\end{numcases}
where $\displaystyle \lambda_0=\frac{\eps}{\mathfrak{Fr}^2}=\mathcal{O}(1)$ and $\displaystyle \gamma_0 = \frac{2\alpha_0(\phi_{\mathrm{max}}-\phi_0)}{\mathfrak{Di}(\phi_0-\mathrm{min})} = \mathcal{O}(1)$.

In this system, the evolution of the volume fraction~$\psi$ is decoupled from the velocity-stress system~\eqref{syst1qua}--\eqref{syst4qua}.
This independent system on~$(\v,q,\btau)$ is similar to the system~\eqref{pb1}--\eqref{pb3} (in which viscosity has been added).
In other words, the system studied in the first sections can be seen as an approximation of the complete system with the $\mu(I)$-rheology.

\subsection{Other models and other possible works}\label{subsec:toto}

\paragraph{Case of a granular material immersed in a liquid}
Note that in a recent paper by Barker et al~\cite{barker23}, the authors use the $\mu(J)$-rheology, instead of $\mu(I)$, to describe fluidised granular flow, i.e. granular material immersed in water.
The dimensionless number~$J$, defined by
$J = \eta_f|\S|/p$ is used for granular flows with a low Stokes number ($St = \rho_0 d^2|\S|/\eta_f$) and is therefore well suited for granular flows in liquids.
As in the case of the models using the $\mu(I)$-rheology discussed in the introduction, instability problems are frequent, and authors often add regularization terms without giving rigorous arguments for the stability of the resulting system. For example, the models introduced in~\cite{montella23,montella21} are relatively close to what is proposed here, but the adapted closure (by adding a time derivative on pressure terms) does not seem to ensure stability or the existence of a solution.
The study proposed in the present article might be adaptable to such configuration.

\paragraph{Shallow water model}
In gravity flow applications, the domain geometry can be highly anisotropic. Typically, pyroclastic flow lengths are of the order of several kilometers, while their heights do not exceed a few meters.
In this context, many authors use models averaged over the transverse component of the flow, thus reducing the number of unknowns and variables in the problem, see for instance~\cite{bouchut15,bouchut16,bouchut21}.
The natural question, then, is what the model introduced in this article might yield when averaged vertically.

\paragraph{Numerical simulations}
The aim of this type of model is to go as far as possible and run numerical simulations to compare with experience.
Actually work is in progress to develop a numerical method, in line with the theoretical results presented here, i.e. retaining as many of the proven properties as possible (bounds on volume fraction, dissipated energy, etc.).


\section*{Acknowledgments}
This is contribution no. xxx of the ClerVolc program of the International Research Center for Disaster Sciences and Sustainable Development of the University of Clermont Auvergne.
\medskip

This project was supported by the French institute of Mathematics for Planet Earth (iMPT) and by the project ComplexFlows of the PEPR Math-VivEs, ANR-23-EXMA- 0004.\medskip

Declaration of Interests: The authors report no conflict of interest.

\bibliography{biblio}
\bibliographystyle{abbrv}

\end{document}